\documentclass[12pt]{amsart}

\usepackage{amsmath}
\usepackage{amssymb}
\usepackage[all]{xy}

\numberwithin{equation}{section}

\newtheorem{thm}{Theorem}[section]

\newtheorem{lem}[thm]{Lemma}
\newtheorem{prop}[thm]{Proposition}
\newtheorem{cor}[thm]{Corollary}
\newtheorem{rem}[thm]{Remark}

\newtheorem{exam-nota}[thm]{Example-Notation}
\newtheorem{nota}[thm]{Notation}
\newtheorem{dfn}[thm]{Definition}

\newtheorem{dfn-nota}[thm]{Definition-Notation}

\newtheorem{dfn-lem}[thm]{Lemma-Definition}

\newcommand{\beqa}{\begin{eqnarray*}}
\newcommand{\eeqa}{\end{eqnarray*}}

\newcommand{\fk}{\mbox{${\mathfrak k}$}}
\newcommand{\fg}{\mbox{${\mathfrak g}$}}
\newcommand{\fq}{\mbox{${\mathfrak q}$}}
\newcommand{\fl}{\mbox{${\mathfrak l}$}}

\newcommand{\fh}{\mbox{${\mathfrak h}$}}
\newcommand{\fn}{\mbox{${\mathfrak n}$}}

\newcommand{\fp}{\mbox{${\mathfrak p}$}}
\newcommand{\fr}{\mbox{${\mathfrak r}$}}
\newcommand{\fb}{\mbox{${\mathfrak b}$}}
\newcommand{\fz}{\mbox{${\mathfrak z}$}}
\newcommand{\fm}{\mbox{${\mathfrak m}$}}
\newcommand{\fu}{\mbox{${\mathfrak u}$}}

\newcommand{\C}{\mbox{${\mathbb C}$}}

\newcommand{\Ad}{{\rm Ad}}

\newcommand{\fgl}{\mathfrak{gl}}

\newcommand{\dn}{{n\choose 2}}

\newcommand{\B}{\mathcal{B}}

\newcommand{\Ybij}{Y_{{\fb}_{i,j}}}
\newcommand{\Ypij}{Y_{{\fp}_{i,j}}}

\newcommand{\bundles}{T^{*}(\B_{n})\times T^{*}(\B_{n-1})}
\newcommand{\X}{\mathfrak{X}}

\title{ Eigenvalue coincidences and $K$-orbits, I}

\author[M. Colarusso]{Mark Colarusso}
\address{Department of Mathematical Sciences, University of Wisconsin-Milwaukee, Milwaukee, WI, 53201}
\email{colaruss@uwm.edu}

\author[S. Evens]{Sam Evens}
\address{Department of Mathematics, University of Notre Dame, Notre Dame, IN, 46556}
\email{sevens@nd.edu}

\begin{document}
\maketitle

\begin{abstract}
We study the variety $\fg(l)$ consisting of matrices $x \in \fgl(n,\C)$ 
such that $x$ and its $n-1$ by $n-1$ cutoff $x_{n-1}$ share exactly
$l$ eigenvalues, counted with multiplicity.   We determine
the irreducible components of $\fg(l)$ by using the
orbits of $GL(n-1,\C)$ on the flag variety $\B$ of $\fgl(n,\C)$.
More precisely, let $\fb \in \B$ be a Borel subalgebra such
that the orbit $GL(n-1,\C)\cdot \fb$ in $\B$ has codimension $l$.
  Then we show that the set $Y_{\fb}:= \{ \Ad(g)(x): x\in \fb \cap
\fg(l), g\in GL(n-1,\C)\}$ is an irreducible component of $\fg(l)$,
and every irreducible component of $\fg(l)$ is of the form
$Y_{\fb}$, where $\fb$ lies in a $GL(n-1,\C)$-orbit of codimension
$l$.   An important ingredient in our proof is the flatness of
a variant of a morphism considered by Kostant and Wallach,
and we prove this flatness assertion using ideas from symplectic
geometry.

\end{abstract}

\section{introduction}\label{sect_intro}
Let $\fg:=\fgl(n,\C)$ be the Lie algebra of $n\times n$ complex matrices.
For $x \in \fg$, let $x_{n-1} \in \fgl(n-1,\C)$ be the upper left-hand
$n-1$ by $n-1$ corner of the matrix $x$.   For $0\le l \le n-1$, we consider
the subset $\fg(l)$ consisting of elements $x\in \fg$ such that
$x$ and $x_{n-1}$ share exactly $l$ eigenvalues, counted with multiplicity.
In this paper, we study the algebraic geometry of the set $\fg(l)$
using the orbits of $GL(n-1,\C) \times GL(1,\C)$ on the flag variety
$\B$ of Borel subalgebras of $\fg$.   In particular, we determine the
irreducible components of $\fg(l)$ and use this to describe
 elements of $\fg(l)$ up to $GL(n-1,\C) \times GL(1,\C)$-conjugacy.

In more detail, let $G=GL(n,\C)$ and let $\theta: G\to G$ be the involution 
$\theta(x)=dxd^{-1}$, where $d=\mbox{diag}[1,\dots, 1,-1]$.  Let $K:=G^{\theta}=GL(n-1,\C)\times GL(1,\C)$.
 It is well-known that $K$ has exactly $n$ closed orbits
on the flag variety $\B$, and each of these closed orbits is isomorphic
to the flag variety $\B_{n-1}$ of Borel subalgebras of $\fgl(n-1,\C)$.
Further, there are finitely many $K$-orbits on $\B$, and for
each of these $K$-orbits $Q$, we consider its length
 $l(Q)=\dim(Q) - \dim(\B_{n-1})$.
It is elementary to verify that $0\le l(Q) \le n-1$.
For $Q=K\cdot \fb_Q$, we consider the $K$-saturation $Y_{Q}:=\Ad(K) \fb_Q$
of $\fb_Q$, which is independent of the choice of $\fb_Q \in Q$.

 \begin{thm}\label{thm:thetheorem}
The irreducible component decomposition of $\fg(l)$ is
 \begin{equation}\label{eq:theorem}
\fg(l) = \bigcup_{l(Q)=n-1-l} Y_Q \cap \fg(l).
\end{equation}
 \end{thm}

The proof uses several ingredients.   The first is the flatness
of a variant of a morphism studied by Kostant and Wallach \cite{KW1}, which implies that  
 $\fg(l)$ is equidimensional.   We prove the flatness assertion using
dimension estimates derived from symplectic geometry, but it also 
follows from results of Ovsienko and Futorny \cite{Ov}, \cite{FuOvsfa}.   
The remaining ingredient is an explicit description of the $l+1$
$K$-orbits $Q$ on $\B$ with $l(Q)=n-1-l$, and the closely related study of
$K$-orbits on generalized flag varieties $G/P$.   Our theorem has the
following consequence.
 Let $\fb_+$ denote the Borel
subalgebra consisting of upper triangular matrices. For $i=1, \dots, n$,
let $(i\,n)$
be the permutation matrix corresponding to the transposition interchanging
$i$ and $n$, and let $\fb_{i}:= \Ad(i\,n)\fb_+$.

\begin{cor}\label{cor:kconjlength}
If $x\in \fg(l)$, then $x$ is $K$-conjugate to an element
in one of $l+1$ explicitly determined $\theta$-stable parabolic subalgebras.   In particular,
if $x \in \fg(n-1)$, then $x$ is $K$-conjugate to an element of $\fb_i$,
where $i=1, \dots, n$.
\end{cor}

This paper is part of a series of papers on $K$-orbits on $\B$ and
the Gelfand-Zeitlin system.  In \cite{CEKorbs}, we used $K$-orbits to
determine the so-called strongly regular elements in the nilfiber of the moment map of 
the Gelfand-Zeitlin system.  These are matrices $x\in\fg$ such that $x_{i}$ is 
nilpotent for all $i=1,\dots, n$ with the added condition that the differentials of the Gelfand-Zeitlin 
functions are linearly independent at $x$.  The strongly regular elements were first studied
extensively in \cite{KW1}.  In later work, we will refine Corollary
\ref{cor:kconjlength} to provide a standard form for all elements of
$\fg(l)$.   This uses $K$-orbits and a finer study of the algebraic
geometry of the varieties $\fg(l)$.  In particular, we will give a more
conceptual proof of the main result from \cite{Col1} and use $K$-orbits 
to describe the geometry of arbitrary fibers of the moment map for 
the Gelfand-Zeitlin system.  

 The work by the  
second author was partially supported by NSA grant 
 H98230-11-1-0151. We would like to thank Adam Boocher and Claudia
Polini for useful discussions.

\section{Preliminaries}\label{sec_prelim}

We show flatness of the partial Kostant-Wallach morphism and recall
needed results concerning $K$-orbits on $\B$.

\subsection{The partial Kostant-Wallach map}
\label{sec_partialkw}

For $x\in \fg$ and $i=1,\dots n$, let $x_i \in \fgl(i,\C)$ denote the upper left
$i\times i$ corner of the matrix $x$.  For any $y\in\fgl(i,\C)$,
let $tr(y)$ denote the trace of $y$.  For $j=1,\dots, i$, 
let $f_{i,j}(x)=tr((x_i)^j)$, which is a homogeneous
function of degree $j$ on $\fg$.
The Gelfand-Zeitlin collection of functions is the set
 $J_{GZ}=\{ f_{i,j}(x) : i=1, \dots, n,
j=1, \dots, i \}$.  The restriction of these functions to any regular
adjoint orbit in $\fg$ produces an integrable system on the orbit \cite{KW1}.  
Let $\chi_{i,j}:\fgl(i,\C)\to \C$ be the function
$\chi_{i,j}(y)=tr(y^{j})$, so that $f_{i,j}(x)=\chi_{i,j}(x_{i})$ and 
$\chi_{i}:=(\chi_{i,1},\dots, \chi_{i,i})$ is the adjoint quotient for $\fgl(i,\C)$.
The Kostant-Wallach map is the morphism given by
\begin{equation} \label{eq:kwdef}
\Phi:\fg \to \C^1 \times \C^2 \times \dots \times \C^n; \, \Phi(x)=(\chi_{1}(x_{1}),\dots, \chi_{n}(x)).
\end{equation}
We will also consider the partial Kostant-Wallach map given by the morphism
\begin{equation} \label{eq:partialkwdef}
\Phi_{n}:\fg\to\C^{n-1}\times\C^{n};\, \Phi_n(x)=(\chi_{n-1}(x_{n-1}), \chi_n(x)).
\end{equation}
Note that
\begin{equation}\label{eq:kwvspartial}
\Phi_n = pr \circ \Phi,
\end{equation}
where $pr:\C^1 \times \C^2 \times \dots \times \C^n \to \C^{n-1} \times \C^n$
is projection on the last two factors.   

\begin{rem}\label{rem_kwsurj}
By Theorem 0.1 of \cite{KW1}, the map $\Phi$ is surjective, and it follows
easily that $\Phi_n$ is surjective.
\end{rem}

We let $I_{n}=(\{f_{ij}\}_{i=n-1, n; j=1, \dots, i})$ denote the ideal generated
by the functions $J_{GZ,n}:=\{f_{i,j}: i=n-1,\, n;\; j=1,\dots, i\}.$  We call the vanishing set 
$V(I_{n})$ the variety of \emph{partially strongly nilpotent matrices} and denote it by $SN_{n}$. 
Thus, 
\begin{equation}\label{eq:strongnil}
SN_{n}:=\{x\in\fg: \, x,\, x_{n-1}\mbox{ are nilpotent}\}.  
\end{equation}
We let $\Gamma_{n}:=\C [\{f_{ij}\}_{i=n-1, n; j=1, \dots, i}]$ be the subring of regular 
functions on $\fg$ generated by $J_{GZ,n}$.

Recall that if $Y \subset \C^m$ is a closed equidimensional subvariety of
dimension
$m-d$, then $Y$ is called a complete intersection if $Y=V(f_1, \dots, f_d)$
is the vanishing set of $d$ functions.   
\begin{thm}\label{thm:completeinter}
The variety of partially strongly nilpotent matrices $SN_{n}$ is a complete intersection of dimension 
\begin{equation}\label{eq:dn}
d_{n}:=n^{2}-2n+1.
\end{equation}
\end{thm}
Before proving Theorem \ref{thm:completeinter}, we show how it implies the flatness 
of the partial Kostant-Wallach map $\Phi_{n}$.

\begin{prop} \label{prop_kwflat}
\begin{enumerate}
\item For all $c\in \C^{n-1}\times \C^{n}$, $\dim(\Phi_{n}^{-1}(c))=n^2-2n+1.$
Thus, $\Phi^{-1}_{n}(c)$ is a complete intersection.  
\item The partial Kostant-Wallach map 
 $\Phi_n:\fg \to \C^{2n-1}$ is a flat morphism.  Thus, $\C[\fg]$ is flat 
 over $\Gamma_{n}$.  
\end{enumerate}
\end{prop}

\begin{proof}
For $x\in \fg$, we let $d_x$ be the maximum of the dimensions of
irreducible components of $\Phi_{n}^{-1}(\Phi_{n}(x))$.   For $c\in \C^{n-1}\times\C^{n}$,
 each irreducible component of
$\Phi_{n}^{-1}(c)$ has dimension at least $d_{n}$ since
$\Phi_{n}^{-1}(c)$ is defined by $2n-1$ equations in $\fg$.
Hence, $d_x \ge d_{n}$.  Since the functions $f_{i,j}$
are homogeneous, it follows that scalar multiplication by
$\lambda \in \C^{\times}$ induces an isomorphism $\Phi_{n}^{-1}(\Phi_{n}(x)) \to
\Phi_{n}^{-1}(\Phi_{n}(\lambda x))$.  
It follows that $d_x = d_{\lambda x}$.
By upper semi-continuity of dimension (see for example, Proposition 4.4 of \cite{Hum}), the set of $y \in \fg$ such that $d_y \ge d$ is closed for
each integer $d$.   It follows that $d_0 \ge d_x$.  
By Theorem \ref{thm:completeinter}, $d_0=d_{n}$.  The first assertion
follows easily.   The second assertion now follows by
 the corollary to Theorem 23.1 of \cite{Mat}. 
\end{proof}

\begin{rem}
\label{r.freeness}
We note that Proposition \ref{prop_kwflat} implies that 
$\C[\fg]$ is free over $\Gamma_n$.  This follows from a
 result in commutative algebra.  Let $A=\oplus_{n\ge 0} A_n$ be
a graded ring
with $A_0=k$ a field
and let $M=\oplus_{n\ge 0} M_n$ be a graded $A$-module.  The needed result
asserts that
$M$ is flat over $A$ if and only if $M$ is free over $A$.  One direction
of this assertion is obvious, and the other direction may be proved
using the same argument as in the proof of Proposition 20 on page
73 of \cite{Serrelocalg}, which is the analogous assertion for
finitely generated modules over local rings.   In \cite{Serrelocalg},
the assumption that $M$ is finitely generated over $A$ is needed
only to apply Nakayama's lemma, but in our graded setting, Nakayama's lemma
(with ideal $I=\oplus_{n > 0} A_n$) does not require the module
$M$ to be finitely generated.
\end{rem}

\begin{rem}
\label{r.fuovref}
Let $I=(J_{GZ})$ be the ideal in $\C[\fg]$ generated by the Gelfand-Zeitlin collection of functions
$J_{GZ}$, and let $SN=V(I)$ be the strongly nilpotent matrices, i.e.,
$SN = \{ x\in \fg : x_i \ \text{is nilpotent for} \ i=1, \dots, n \}$.
Ovsienko proves in \cite{Ov} that $SN$ is a complete intersection,
and results of Futorny and Ovsienko from \cite{FuOvsfa} show that
Ovsienko's theorem implies that
$\C[\fg]$ is free over $\Gamma:=\C[\{f_{i,j}\}_{i=1, \dots, n;\, j=1, \dots, i}]$.
It then follows easily that $\C[\fg]$ is flat over $\Gamma_n$, and hence
that $\Phi_n$ is flat.  Although we could have simply cited the
results of Futorny and Ovsienko to prove flatness of $\Phi_n$, we
prefer our approach, which we regard as more conceptual.
\end{rem}

\begin{proof} [Proof of Theorem \ref{thm:completeinter}]
Let $\X$ be an irreducible component of $SN_{n}$.  
We observed in the proof of Proposition \ref{prop_kwflat} that $\dim \X\geq d_{n}$.  To show $\dim \X\leq d_{n}$, we consider a generalization of the Steinberg variety (see Section 3.3 of \cite{CG}).  We first recall a few facts about the cotangent bundle to the flag variety.

For the purposes of this proof, we denote the flag variety of $\fgl(n,\C)$ by $\B_{n}$.  
We consider the form $<<\cdot, \cdot>>$ on $\fg$ given by 
$<<x,y>>=tr(xy)$ for $x,\, y\in\fg$.  If $\fb\in\B_{n}$, the annihilator $\fb^{\perp}$ of 
$\fb$ with respect to the form $<<\cdot, \cdot>>$ is $\fn=[\fb,\fb]$.  
We can then identify $T^{*}(\B_{n})$ with the closed subset of 
$\fg\times \B_{n}$ given by:
$$
T^{*}(\B_{n})=\{(x,\fb): \, \fb\in \B_{n}, x\in\fn\}.  
$$

We let $\fg_{n-1}=\fgl(n-1,\C)$ and view $\fg_{n-1}$ as a subalgebra of $\fg$ by embedding $\fg_{n-1}$ in the top lefthand 
corner of $\fg$.  Since $\fg$ is the direct sum $\fg=\fg_{n-1}\oplus\fg_{n-1}^{\perp}$, the restriction of 
$<<\cdot, \cdot>>$ to $\fg_{n-1}$ is non-degenerate.  For a Borel subalgebra $\fb^{\prime} \in \B_{n-1}$, we let $\fn^{\prime}=[\fb^{\prime}, \fb^{\prime}]$.  
We consider a closed subvariety $Z \subset \fg \times \B_n \times \B_{n-1}$ 
defined as follows: 
\begin{equation}\label{eq:newstein}
Z=\{(x,\fb, \fb^{\prime}):  \fb\in \B_{n}\; , \, \fb^{\prime}\in\B_{n-1}\mbox { and } x\in \fn,\, x_{n-1}\in \fn^{\prime}\}.
\end{equation}
Consider the morphism $\mu: Z\to \fg$, where $\mu(x,\fb,\fb^{\prime})=x.$ 
Since the varieties $\B_{n}$ and $\B_{n-1}$ are projective, the morphism $\mu$ is proper.  

We consider the closed embedding $Z \hookrightarrow \bundles \cong T^{*}(\B_{n}\times \B_{n-1})$ given by
$(x, \fb, \fb^{\prime}) \to (x,-x_{n-1}, \fb, \fb^{\prime})$.  We denote the image of $Z$ under this embedding
by $\widetilde{Z}\subset T^{*}(\B_{n}\times \B_{n-1})$.  Let $G_{n-1}$ be the closed subgroup of $GL(n,\C)$  corresponding to $\fg_{n-1}$. 
Then $G_{n-1}$ acts diagonally on $\B_{n}\times \B_{n-1}$ via $k\cdot(\fb,\fb^{\prime})=(k\cdot\fb,k\cdot\fb^{\prime})$ for $k\in G_{n-1}$.   We claim $\widetilde{Z}\subset T^{*}(\B_{n}\times \B_{n-1})$ is the union of conormal bundles to the $G_{n-1}$-diagonal orbits 
in $\B_{n}\times \B_{n-1}$.  Indeed, let $(\fb,\fb^{\prime})\in\B_{n}\times\B_{n-1}$, and let $Q$ be its $G_{n-1}$-orbit.  
Then 
$$
T_{\left(\fb,\fb^{\prime}\right)}(Q)=\mbox{span}\{(Y\;\mbox{mod } \fb, \, Y\;\mbox{mod } \fb^{\prime}):\, Y\in\fg_{n-1}\}.  
$$
Now let $(\lambda_{1},\lambda_{2})\in (\fn,\fn^{\prime})$ with $(\lambda_{1},\lambda_{2})\in (T^{*}_{Q})(\B_{n}\times \B_{n-1})_{\left(\fb,\fb^{\prime}\right)}$, the
fiber of the conormal bundle to $Q$ in $\B_{n} \times \B_{n-1}$ at
the point $(\fb, \fb^{\prime})$.
Then 
$$
<<\lambda_{1},Y>>+<<\lambda_{2},Y>>=0 \mbox{ for all } Y\in\fg_{n-1}.
$$
Thus, $\lambda_{1}+\lambda_{2}\in\fg_{n-1}^{\perp}$.  But since 
$\lambda_{2}\in\fn^{\prime}\subset\fg_{n-1}$, it follows that $\lambda_{2}=-(\lambda_{1})_{n-1}$. 
Thus, 
$$
T^{*}_{Q}(\B_{n}\times \B_{n-1})=\{(\mu_{1},\fb_{1},-(\mu_{1})_{n-1},\fb_{2}),\; \mu_{1}\in \fn_{1},\, (\mu_{1})_{n-1}\in\fn_{2},\mbox{ where }
(\fb_{1},\fb_{2})\in Q\}.
$$

We recall the well-known fact that there are only
 finitely many $G_{n-1}$-diagonal orbits in $\B_{n}\times \B_{n-1}$,
which follows from \cite{KimelVin}, \cite{Brionclass}, 
or in a more explicit form is proved in \cite{Hashi}.
Therefore, the irreducible component decomposition of $\widetilde{Z}$ is:
$$
\widetilde{Z}=\bigcup_{i} \overline{T^{*}_{Q_{i}}(\B_{n}\times\B_{n-1})}\subset T^{*}(\B_{n}\times \B_{n-1}),
$$ 
where $i$ runs over the distinct $G_{n-1}$-diagonal orbits in $\B_{n}\times \B_{n-1}$.  Thus, $Z\cong \widetilde{Z}$ is a closed, equidimensional subvariety of dimension 
$$
\dim Z=\frac{1}{2} (\dim T^{*}(\B_{n}\times \B_{n-1}) )=d_{n}.
$$

Note that $\mu:Z\to SN_{n}$ is surjective. 
Since $\mu$ is proper, for every irreducible component 
$\X\subset SN_{n}$ of $SN_{n}$,
we see that
\begin{equation}\label{eq:acptofSNn}
\X=\mu(Z_{i})
\end{equation}
for some irreducible component $Z_{i}\subset Z$.
Since $\dim Z_{i}=d_{n}$ and $\dim \X\geq d_{n}$, we conclude that $\dim \X=d_{n}$. 
\end{proof}

In Proposition \ref{prop:irredsnn}, we will determine the irreducible
components of $SN_{n}$ explicitly.

\subsection{$K$-orbits}
\label{sec_korbitsprelim}

We recall some basic facts about $K$-orbits  on generalized flag
varieties $G/P$ (see \cite{Mat79, RS, MO, Yam, CEexp} for more details).

By the general theory of orbits of symmetric subgroups on generalized
flag varieties, $K$ has finitely many orbits on
$\B$.   For this paper, it is useful to parametrize the orbits.
To do this, we let $B_+$ be the upper triangular Borel subgroup of $G$,
and identify $\B \cong G/B_+$ with the variety of flags in $\C^n$.
 We use the following notation for flags in $\C^{n}$.  Let 
   $$
  \mathcal{F}=( V_{0}=\{0\}\subset V_{1}\subset\dots\subset V_{i}\subset\dots\subset V_{n}=\C^{n}).
   $$
   be a flag in $\C^{n}$, with $\dim V_{i}=i$ and $V_{i}=\mbox{span}\{v_{1},\dots, v_{i}\}$, with each $v_{j}\in\C^{n}$.  We will denote the flag $\mathcal{F}$ as follows:
   $$
   v_{1}\subset  v_{2}\subset\dots\subset v_{i}\subset v_{i+1}\subset\dots\subset v_{n}. 
   $$
    We denote the standard ordered basis of $\C^{n}$ by $\{e_{1},\dots, e_{n}\}$,
and let $E_{i,j} \in \fg$ be the matrix with $1$
in the $(i,j)$-entry and $0$ elsewhere.  
  
   

 There are 
$n$ closed $K$-orbits on $\B$ (see Example 4.16 of \cite{CEexp}),
 $Q_{i,i}=K\cdot \fb_{i,i}$ for $i=1,\dots, n$,  where 
the Borel subalgebra $\fb_{i,i}$ 
 is the stabilizer of the following flag in $\C^{n}$:
   \begin{equation}\label{eq:flagi}
 \mathcal{F}_{i,i}=( e_{1}\subset \dots\subset e_{i-1}\subset \underbrace{e_{n}}_{i}\subset e_{i}\subset \dots\subset e_{n-1}).
   \end{equation}
     Note that if $i=n$, then the flag $\mathcal{F}_{i,i}$ is the standard flag $\mathcal{F}_{+}$:
\begin{equation} \label{eq:stdflag}
   {\mathcal{F}}_+ = (e_{1} \subset \dots \subset e_{n}),
\end{equation}
and $\fb_{n,n}=\fb_{+}$ is the standard Borel subalgebra of $n\times n$ upper triangular matrices. It is not difficult to check that 
$K\cdot \fb_{i,i} = K\cdot \Ad(i\,n) \fb_+$.  
If $i=1$, then $K\cdot \fb_{1,1}=K\cdot \fb_{-}$, where $\fb_{-}$ is the Borel subalgebra of lower triangular matrices in $\fg$.

The non-closed $K$-orbits in $\B$ are the orbits
$Q_{i,j}=K\cdot \fb_{i,j}$
 for $1\leq i<j\leq n$, where $\fb_{i,j}$ 
is the stabilizer of the flag in $\C^{n}$:
  \begin{equation} \label{eq:flagij}
  \mathcal{F}_{i,j} =   (e_{1}\subset \dots\subset \underbrace{e_{i}+e_{n}}_{i}\subset e_{i+1}\subset \dots \subset e_{j-1} \subset \underbrace{e_{i}}_{j}\subset e_{j} \subset \dots \subset e_{n-1}).
\end{equation}
 There are $\dn$ such orbits (see Notation
4.23 and Example 4.31 of \cite{CEexp}).
 
 Let  $w$ and $\sigma$ be the permutation matrices 
corresponding respectively to the cycles  $(n\,n-1\cdots\, i)$ and $ (i+1\, i+2\dots j)$,
and let $u_{\alpha_{i}}$ be the Cayley transform matrix such that
$$
u_{\alpha_{i}}(e_{i})=e_{i} + e_{i+1}, \ 
u_{\alpha_{i}}(e_{i+1})=-e_{i} + e_{i+1}, \
u_{\alpha_{i}}(e_{k})=e_{k}, k\not= i, i+1.
$$
 For $1\leq i\leq j\leq n$, we define: 
\begin{equation} \label{eq:vdef} 
v_{i,j}:=\left\{ \begin{array}{ccc} w& \mbox{ if }& i=j\\
wu_{\alpha_{i}}\sigma& \mbox{ if } &i\neq j\end{array}\right\}
\end{equation} 
It is easy to verify that $v_{i,j}({\mathcal{F}}_+)={\mathcal{F}}_{i,j}$, and thus
$\Ad(v_{i,j})\fb_{+}=\fb_{i,j}$ (see Example 4.30 of \cite{CEexp}).

 \begin{rem}\label{r:length}
 The length of the $K$-orbit $Q_{i,j}$ 
 is $l(Q_{i,j})=j-i$ for any $1\leq i\leq j\leq n$ (see Example 4.30 of \cite{CEexp}).  For example, a $K$-orbit $Q_{i,j}$ is closed 
 if and only if $Q=Q_{i,i}$ for some $i$.   The $n-l$ orbits of length $l$ are $Q_{i,i+l}$,
$i=1, \dots, n-l$. 
\end{rem}

For a parabolic subgroup $P$ of $G$ with Lie algebra $\fp$, we consider the
generalized flag variety $G/P$, which we identify with parabolic subalgebras
 of type
$\fp$ and with partial flags of type $\fp$.   
 We will make use of the following notation for partial flags.  Let
   $$
   \mathcal{P}=(V_{0}=\{0\}\subset V_{1}\subset\dots\subset V_{i}\subset\dots\subset V_{k}=\C^{n})
   $$
   denote a $k$-step partial flag with $\dim V_{j}=i_{j}$ and $V_{j}=\mbox{span}\{v_{1},\dots, v_{i_{j}}\}$ for $j=1,\dots, k$.  
   Then we denote $\mathcal{P}$ as 
   $$
   v_{1},\dots, v_{i_{1}} \subset v_{i_{1}+1},\dots, v_{i_{2}}\subset\dots\subset v_{i_{k-1}+1},\dots, v_{i_{k}}.
   $$

  In particular for $i\leq j$, we let $\fr_{i,j}\subset\fg$ denote the parabolic subalgebra which is the stabilizer of the $n-(j-i)$-step partial flag in $\C^{n}$
 \begin{equation}\label{eq:standardpflag}
 \mathcal{R}_{i,j}=(e_{1}\subset e_{2}\subset\dots\subset e_{i-1}\subset e_{i},\dots, e_{j}\subset e_{j+1}\subset\dots\subset e_{n}). 
 \end{equation}
It is easy to see that $\fr_{i,j}$ is the standard parabolic subalgebra generated by the Borel subalgebra $\fb_{+}$ and
the negative simple root spaces 
$\fg_{-\alpha_{i}}, \fg_{-\alpha_{i+1}}, \dots, \fg_{-\alpha_{j-1}}$.
We note that $\fr_{i,j}$ has Levi decomposition $\fr_{i,j} = \fm + \fn$, with $\fm$ consisting
of block diagonal matrices of the form
\begin{equation} \label{eq:levifactor}
\fm=\underbrace{\fgl(1,\C)\oplus\dots\oplus \fgl(1,\C)}_{i-1 \mbox{ factors}}\oplus\fgl(j+1-i,\C)\oplus\underbrace{ \fgl(1,\C)\oplus\dots\oplus \fgl(1,\C)}_{n-j\mbox{ factors}}.
\end{equation}
 Let $R_{i,j}$ be the parabolic subgroup of $G$ with Lie algebra $\mathfrak{r}_{i,j}$.
 Let $\fp_{i,j}:=\Ad(v_{i,j})\fr_{i,j}\in G/R_{i,j}$, where $v_{i,j}$ is defined in (\ref{eq:vdef}).  Then $\fp_{i,j}$ is the stabilizer of the partial flag
\begin{equation}\label{eq:firstpartial}
\mathcal{P}_{i,j}=(e_{1}\subset e_{2}\subset\dots\subset e_{i-1}\subset e_{i},\dots, e_{j-1}, e_{n}\subset e_{j}\subset\dots\subset e_{n-1}), 
\end{equation}
and $\fp_{i,j}\in G/R_{i,j}$ is a $\theta$-stable parabolic subalgebra of $\fg$.  Indeed, recall that $\theta$ is given by conjugation 
by the diagonal matrix $d=diag[1,\dots, 1,-1]$.  Clearly $d(\mathcal{P}_{i,j})=\mathcal{P}_{i,j}$, whence $\fp_{i,j}$ is $\theta$-stable. Moreover, the parabolic subalgebra $\fp_{i,j}$ has Levi decomposition 
$\fp_{i,j}=\fl\oplus\fu$ where both $\fl$ and $\fu$ are $\theta$-stable and $\fl$ is isomorphic to the Levi subalgebra in Equation (\ref{eq:levifactor}).
Since $\fp_{i,j}$ is $\theta$-stable, it follows from Theorem 2 of \cite{BH} that the $K$-orbit  $Q_{\fp_{i,j}}=K\cdot \fp_{i,j}$ is closed in  $G/R_{i,j}$.  

For a parabolic subgroup $P\subset G$ with Lie algebra $\fp\subset\fg$,
consider the partial Grothendieck resolution
$\tilde{\fg}^{\fp}=\{ (x,\fr)\in \fg \times G/P \; | \; x\in\fr\}$,
as well as the morphisms  
$\mu:\tilde{\fg}^{\fp}\to \fg,\; \mu(x,\fr)=x$, and $\pi: \tilde{\fg}^{\fp}\to G/P,\; \pi(x,\fr)=\fr$.  
Then $\pi$ is a smooth morphism of relative dimension $\dim \fp$ 
(for $G/B$, see Section 3.1 of \cite{CG} and Proposition III.10.4 of \cite{Ha}, and
the general case of $G/P$ follows by the same argument). 
   For $\fr \in G/P$, let $Q_{\fr}=K\cdot \fr \subset G/P$.
Then $\pi^{-1}(Q_{\fr})$ has dimension $\dim(Q_{\fr}) + \dim(\fr)$.
It is well-known that $\mu$ is proper and its restriction to $\pi^{-1}(Q_{\fr})$ 
generically has finite fibers (Proposition 3.1.34 and Example
3.1.35 of \cite{CG} for the case of $G/B$, and again the general case has a similar
proof).  

\begin{nota} For a parabolic subalgebra $\fr$ with $K$-orbit $Q_{\fr}\in G/P$, we consider the irreducible
subset
\begin{equation}\label{eq:Yr}
Y_{\fr}:=\mu(\pi^{-1}(Q_{\fr}))=\Ad(K)\fr.
\end{equation}
To emphasize the orbit $Q_{\fr}$, we will also denote this set as
\begin{equation}\label{eq:YQps}
Y_{Q_{\fr}}:=Y_{\fr}.
\end{equation}
\end{nota}

 It follows from generic finiteness of $\mu$ that $Y_{Q_{\fr}}$
contains an open subset of dimension
\begin{equation}\label{eq:YQdim}
\dim(Y_{Q_{\fr}}):=\dim \pi^{-1}(Q_{\fr})=\dim \fr+\dim(Q_{\fr})=\dim\fr+\dim (\fk/ \fk\cap \fr),
\end{equation}
where $\fk=Lie(K)=\fgl(n-1,\C)\oplus\fgl(1,\C)$.

\begin{rem} \label{r:closedorbitcase}
Since $\mu$ is proper, when $Q_{\fr}=K\cdot \fr$ is closed in $G/P$, then
$Y_{Q_{\fr}}$ is closed.
\end{rem}
\begin{rem}\label{r:Qstrat}
Note that
$$
\fg=\bigcup_{Q\subset G/P} Y_{Q},
$$
is a partition of $\fg$, where the union 
is taken over the finitely many $K$-orbits in $G/P$. 
\end{rem}

\begin{lem} \label{l:YQclosure}
Let $Q \subset G/P$ be a $K$-orbit.   Then
\begin{equation}
 \overline{Y_{Q}}=\bigcup_{Q^{\prime}\subset \overline{Q}} Y_{Q^{\prime}}.
\end{equation}
\end{lem}

\begin{proof}
Since $\pi$ is a smooth morphism, it is flat by Theorem III.10.2 of \cite{Ha}.  Thus,
by Theorem VIII.4.1 of \cite{SGA1}, 
$\pi^{-1}(\overline{Q})=\overline{\pi^{-1}(Q)}$.   The result follows
since $\mu$ is proper.
\end{proof}

\subsection{Comparison of $K\cdot \fb_{i,j}$ and $K\cdot \fp_{i,j}$}

We prove a technical result that will be needed to prove our main theorem.

\begin{rem} \label{rem:fbijfpij}
Note that $\fb_{i,j} \subset \fp_{i,j}$ and when $i=j$, $\fp_{i,i}$
is the Borel subalgebra $\fb_{i,i}$.   To check the first assertion,
note that $\fb_{+}\subset\fr_{i,j}$ so that $\fb_{i,j}=\Ad(v_{i,j})\fb_{+}\subset \Ad(v_{i,j})\fr_{i,j}=\fp_{i,j}$.  
The second assertion is verified by noting that when $i=j$,
the partial flag $\mathcal{P}_{i,j}$ is the full flag $\mathcal{F}_{i,i}$.
\end{rem}

\begin{prop}\label{p:diffYQs}
Consider the $K$-orbits $Q_{i,j} =K\cdot \fb_{i,j} \subset \B$ and $Q_{\fp_{i,j}}=K\cdot \fp_{i,j} \subset
G/P_{i,j}$, with $1\leq i \leq j \leq n$.   
Then $\dim(\Ybij)=\dim(\Ypij)$ and $\overline{\Ybij}=\Ypij$. 
\end{prop} 

\begin{proof}
By definitions and Remark \ref{rem:fbijfpij}, $\Ybij$ is a constructible subset
of $\Ypij$.   Since $\Ypij$ is closed by Remark \ref{r:closedorbitcase},
and irreducible by construction, it suffices to show that
$\dim(\Ybij)=\dim(\Ypij)$.  

We compute the dimension 
of $\Ybij$ using Equation (\ref{eq:YQdim}).  Since $l(Q_{i,j})=j-i$, it follows
that $\dim Q_{i,j}=\dim\mathcal{B}_{n-1}+j-i$.  
Since $\dim(\B_{n-1})={n-1\choose 2}$, Equation 
(\ref{eq:YQdim}) then implies: 
 \begin{equation}\label{eq:largeYQdim}
  \begin{split}
  \dim\Ybij=\dim\fb_{i,j}+\dim\mathcal{B}_{n-1}+l(Q_{i,j})&={n+1\choose 2}+{n-1\choose 2} + l(Q_{i,j})\\
  &=n^2-n+1+j-i.\\
  \end{split}
  \end{equation}

We now compute the dimension of $\Ypij$.  By Equation (\ref{eq:YQdim}),
it follows that 
\begin{equation}\label{eq:secondpardim}
\dim \Ypij=\dim\fp_{i,j}+\dim \fk-\dim (\fk\cap \fp_{i,j}).  
\end{equation}  
Since both $\fl$ and $\fu$ are $\theta$-stable, it follows that 
$\dim \fk\cap\fp_{i,j}=\dim \fk\cap\fl+\dim \fk\cap \fu.$
To compute these dimensions, it is convenient to use the following 
explicit matrix description of the parabolic subalgebra $\fp_{i,j}$,
which follows from Equation (\ref{eq:firstpartial}).

\begin{equation}\label{eq:bigmatrix}
\fp_{i,j}=\left[\begin{array}{ccccccccccc}
a_{11} &\dots &\dots &a_{1i-1} &a_{1i} &\dots &a_{1j-1} & \dots&\dots &a_{1n-1} & a_{1n}\\
0 & \ddots &  & \vdots &\vdots & *& \vdots&*  &* &\vdots & \vdots\\
\vdots& & & a_{i-1i-1} &\vdots &*&  \vdots& * &* &a_{i-1n-1} & a_{i-1n}\\
& & &0 & a_{ii} &\dots &a_{ij-1} & * & *&a_{in-1} &a_{in}\\
&  & &\vdots & \vdots&\ddots &\vdots &* & *& \vdots& \vdots\\
& & & \vdots&a_{ij-1}&\dots & a_{j-1j-1} &\dots &\dots &a_{j-1n-1} & a_{j-1n}\\
& & & 0&0 &\dots& 0 &a_{jj} & \dots&a_{jn-1} &0 \\
& & & \vdots& \vdots &  & \vdots&0  & \ddots &\vdots &\vdots\\
\vdots & & & \vdots & 0 & & 0 & 0 &0  &a_{n-1n-1} & 0  \\
0&\dots & \dots & 0& a_{ni}& \dots&a_{nj-1} & a_{nj}&\dots & a_{nn-1} &a_{nn}\end{array}\right]_{\mbox{.}}
\end{equation}

Using (\ref{eq:bigmatrix}), we observe that $\fk\cap\fl\cong \fgl(1,\C)^{n-j+i}\oplus \fgl(j-i, \C)$, 
so that $\dim \fk\cap\fl=n-j+i+(j-i)^2$.  
Now $\fu=\fu\cap\fk\oplus\fu^{-\theta}$, where 
$\fu^{-\theta}:= \{ x\in \fu : \theta(x)=-x \}$.  
Using (\ref{eq:bigmatrix}),  we see that
$\fu^{-\theta}$ has basis $\{ E_{n,j}, \dots, E_{n, n-1}, E_{1,n}, \dots, E_{i-1, n}\}$,
so  
$\dim\fu^{-\theta}=n-j+i-1$.  Thus, $\dim \fu\cap\fk=\dim\fu-(n-j+i-1)$.  Putting these observations together, we obtain
\begin{equation}\label{eq:pijdim}
\dim \fk\cap\fp_{i,j}=(j-i)^2+\dim\fu+1.
\end{equation}

Now 
$$\dim\fp_{i,j}=\dim\fl+\dim\fu=(j-i+1)^2+n-j+i-1+\dim\fu.$$
(see Equation (\ref{eq:levifactor})).  Thus, Equation (\ref{eq:secondpardim}) 
implies that 
$$\dim\Ypij=\dim \fk+(j-i+1)^2+n-j+i-1-(j-i)^2-1=n^2-n+1+j-i,$$
which agrees with (\ref{eq:largeYQdim}), and hence completes the proof.

\end{proof}
 
\begin{rem}\label{r:projpij}
It follows from Equation (\ref{eq:bigmatrix}) that
$(\fp_{i,j})_{n-1}:=\pi_{n-1,n}(\fp_{i,j})$ is a parabolic subalgebra, where
$\pi_{n,n-1}:\fg \to \fgl(n-1,\C)$ is the projection
$x\mapsto x_{n-1}$.  Further, with $l=j-i$, 
$(\fp_{i,j})_{n-1}$ has Levi decomposition
$(\fp_{i,j})_{n-1}=\fl_{n-1}\oplus\fu_{n-1}$ with 
$\fl_{n-1}= \fgl(1,\C)^{n-1-l}\oplus \fgl(l,\C)$.
\end{rem}

\section{The varieties $\fg(l)$} \label{sec_g(l)}

In this section, we prove our main results.

For $x\in \fg$, let 
 $\sigma(x)=\{\lambda_{1},\dots, \lambda_{n}\}$ denote its eigenvalues, where an eigenvalue $\lambda$ is listed $k$
times if it appears with multiplicity $k$.
Similarly, let $\sigma(x_{n-1})=\{\mu_{1},\dots, \mu_{n-1}\}$ be the eigenvalues of $x_{n-1} \in \fgl(n-1, \C)$, again listed with multiplicity.
 For $i=n-1, n$, let $\fh_{i}\subset\fg_{i}:= \fgl(i, \C)$ be the standard Cartan subalgebra 
of diagonal matrices.  
We denote elements of $\fh_{n-1}\times\fh_{n}$ by $(x,y)$, with 
$x=(x_{1},\dots, x_{n-1})\in \C^{n-1}$ and $y=(y_{1},\dots, y_{n})\in \C^n$
the diagonal coordinates of $x$ and $y$.
  For $l=0, \dots, n-1$, we define 
\begin{equation*}
\begin{split}
(\fh_{n-1}\times\fh_{n})(\geq l)=&\{(x, y): \exists\; 1\leq i_{1}<\dots <i_{l}\leq n-1
\mbox{ with } x_{i_{j}}=y_{k_{j}}\\
&\mbox{ for some } 1\leq k_{1},\dots, k_{l}\leq n \mbox{ with } k_{j}\neq k_{m}\}.
\end{split}
\end{equation*}
Thus, $(\fh_{n-1}\times\fh_{n})(\geq l)$ consists of elements of $\fh_{n-1}\times\fh_{n}$ with at least $l$ coincidences in the spectrum of $x$ and $y$ counting repetitions.  Note that $(\fh_{n-1}\times\fh_{n})(\geq l)$ is a closed subvariety
of $\fh_{n-1}\times\fh_{n}$ and is equidimensional of codimension
$l$.

 Let $W_{i}=W_{i}(\fg_{i},\fh_{i})$ be the Weyl group of $\fg_{i}$.  Then $W_{n-1}\times W_{n}$ acts on $(\fh_{n-1}\times\fh_{n})(\geq l)$. 
 Consider the finite morphism $p: \fh_{n-1}\times\fh_{n}\to  (\fh_{n-1}\times\fh_{n})/ (W_{n-1}\times W_{n})$.  
Let $F_{i}: \fh_{i}/W_{i} \to \C^{i}$ be the Chevalley isomorphism, and let 
$$V^{n-1,n}:=\C^{n-1}\times\C^{n}, \mbox{ so that } F_{n-1}\times F_{n}: (\fh_{n-1}\times\fh_{n})/ (W_{n-1}\times W_{n}) \to V^{n-1,n}$$ is an isomorphism.  The following varieties play a major role in our study of eigenvalue coincidences.

\begin{dfn-nota}\label{ref:thenotation}
For $l=0,\dots, n-1$, we let 
\begin{equation}\label{eq:Vgel}
V^{n-1,n}(\geq l):=(F_{n-1}\times F_{n})((\fh_{n-1}\times\fh_{n})(\geq l)/(W_{n-1}\times W_{n})),
\end{equation} 
 \begin{equation}\label{eq:Vl}
  V^{n-1,n}(l):=V^{n-1,n}(\geq l)\setminus V^{n-1,n}(\geq l+1).
  \end{equation}
For convenience, we let $V^{n-1,n}(n)=\emptyset$.
\end{dfn-nota}

\begin{lem}\label{lem:vnirreducible}
The set $V^{n-1,n}(\geq l)$ is an irreducible closed subvariety of $V^{n-1,n}$ of
dimension $2n-1-l$.   Further, $V^{n-1,n}(l)$ is open and dense 
in $V^{n-1,n}(\geq l)$.
\end{lem}

\begin{proof}
Indeed, the set $Y:= \{(x,y)\in \fh_{n-1}\times\fh_{n} : x_{i}=y_{i} \mbox{ for }
i=1, \dots, l\}$ is closed and irreducible of dimension $2n-1-l$.   The
first assertion follows 
 since $(F_{n-1} \times F_{n})\circ p$ is a finite morphism and 
$(F_{n-1} \times F_{n})\circ p (Y)=V^{n-1, n}(\geq l)$.
The last assertion of the lemma now follows from Equation (\ref{eq:Vl}).
\end{proof}

\begin{dfn}\label{dfn:ggel}
We let
$$\fg(\geq l):=\Phi_{n}^{-1}(V^{n-1,n}(\geq l)).$$
\end{dfn}

\begin{rem}\label{rem:quotienteigenvalue}
Recall that the quotient morphism 
$p_{i}:\fg_{i} \to \fg_{i}//GL(i,\C) \cong\fh_{i}/W_{i}$
associates to $y \in \fg_{i}$ its spectrum $\sigma(y)$, and 
$(F_{n-1} \times F_{n})\circ (p_{n-1} \times p_{n})=\Phi_{n}$.
It follows that $\fg(\geq l)$ consists of elements of
$x$ with at least $l$ coincidences in the spectrum of
$x$ and $x_{n-1}$, counted with multiplicity.
\end{rem}

 It is routine to check that
 \begin{equation}\label{eq:fgl}
  \fg(l):=\fg(\geq l) \setminus \fg(\geq l+1) = \Phi_n^{-1}(V^{n-1,n}(l))
  \end{equation}
consists of elements of $\fg$ with exactly $l$ coincidences in the
spectrum of $x$ and $x_{n-1}$, counted with multiplicity.

\begin{prop}\label{prop_dimgl}
\begin{enumerate}
\item The variety $\fg(\ge l)$ is equidimensional of dimension $n^2 - l$.
\item $\fg(\ge l) = \overline{\fg(l)} = \bigcup_{k\geq l} \fg(k)$.
\end{enumerate}
\end{prop}

\begin{proof}
By Proposition \ref{prop_kwflat}, the morphism $\Phi_n$ is flat.
By Proposition III.9.5 and Corollary III.9.6 of \cite{Ha}, the variety
$\fg(\geq l)$ is equidimensional of dimension 
$\dim (V^{n-1,n}(\geq l))+(n-1)^{2}$, which
gives the first assertion by Lemma \ref{lem:vnirreducible}.
 For the second assertion, 
by the flatness of $\Phi_n$, Theorem VIII.4.1 of \cite{SGA1},
and Lemma \ref{lem:vnirreducible},
\begin{equation}\label{eq:glbar}
\overline{\fg(l)}=\overline{\Phi_{n}^{-1}(V^{n-1,n}(l))}=\Phi_{n}^{-1}(\overline{V^{n-1,n}(l)})=\Phi_{n}^{-1}(V^{n-1,n}(\geq l))=\fg(\geq l).                        
 \end{equation}
The remaining equality follows from definitions.
\end{proof}


We now relate the partitions $\fg=\bigcup \fg(l)$ and $\fg = \bigcup_{Q\subset\B} Y_Q$ (see Remark \ref{r:Qstrat}). 

\begin{thm}\label{thm:overlaps}
\begin{enumerate}
\item Consider the closed subvarieties $\Ypij$ for $1\leq i\leq j\leq n$, and let $l=j-i$. 
Then $\Ypij\subset\fg(\geq n-1-l)$.
\item  In particular, if $Q\subset\B$ is a $K$-orbit with 
$l(Q)=l$, then $Y_{Q}\subset\fg(\geq n-1-l)$.  
\end{enumerate}
\end{thm}

\begin{proof}
The second statement of the theorem follows from the first statement using
Remark \ref{r:length} and Proposition \ref{p:diffYQs}.

To prove the first statement of the theorem, let $\fq$ be a 
parabolic subalgebra of $\fg$ with $\fq\in Q_{\fp_{i,j}}$, and
let $y\in\fq$.   We need to show that $\Phi_{n}(y)\in V^{n-1,n}(\geq n-1-l)$.  Since the 
map $\Phi_{n}$ is $K$-invariant, it is enough to show that $\Phi_{n}(x)\in V^{n-1,n}(\geq n-1-l)$ 
for $x\in\fp_{i,j}$.  

We recall that $\Phi_{n}(x)=(\chi_{n-1}(x_{n-1}), \chi_{n}(x))$ where $\chi_{i}:\fgl(i,\C)\to \C^{i}$ is the 
adjoint quotient for $i=n-1, n$.  For $x\in\fp_{i,j}$, let $x_{\fl}$ be the projection of 
$x$ onto $\fl$ off of $\fu$.  It is well-known that $\chi_{n}(x)=\chi_{n}(x_{\fl})$.  Using the identification $\fl\cong\fgl(1,\C)^{n-1-l}\oplus \fgl(l+1,\C)$, 
we decompose $x_{\fl}$ as $x_{\fl}=x_{\fgl(1)^{n-1-l}}+x_{\fgl(l+1)}$, where $x_{\fgl(1)^{n-1-l}}\in\fgl(1,\C)^{n-1-l}$ and $x_{\fgl(l+1)}\in\fgl(l+1,\C)$.   It follows that the
coordinates of $x_{\fgl(1)^{n-1-l}}$ are in the spectrum of $x$ (see (\ref{eq:bigmatrix})).

 Recall the 
projection $\pi_{n,n-1}:\fg\to\fg_{n-1}$, $\pi_{n,n-1}(x)=x_{n-1}$.  
Recall the Levi decomposition $(\fp_{i,j})_{n-1}=\fl_{n-1}\oplus\fu_{n-1}$ of
the parabolic subalgebra $(\fp_{i,j})_{n-1}$ of $\fgl(n-1,\C)$ from
Remark \ref{r:projpij}, and recall that
 $\fl_{n-1}= \fgl(1,\C)^{n-1-l}\oplus \fgl(l,\C)$.
Thus,  $\chi_{n-1}(x_{n-1})=\chi_{n-1}((x_{n-1})_{\fl_{n-1}})$.
 We use the decomposition 
$(x_{n-1})_{\fl_{n-1}}=x_{\fgl(1)^{n-1-l}}+\pi_{l+1,l}(x_{\fgl(l+1)})$, where
 $\pi_{l+1, l}:\fgl(l+1,\C)\to \fgl(l,\C)$ 
is the usual projection.  
It now follows easily from Remark \ref{rem:quotienteigenvalue}
that 
$\Phi_n(x)\in V^{n-1,n}(\geq n-1-l)$, since the coordinates of
$x_{\fgl(1)^{n-1-l}}$ are eigenvalues both for $x$ and $x_{n-1}$.
\end{proof}
  

We now recall and prove our main theorem.

 \begin{thm}\label{thm:king}
 Consider the locally closed subvariety $\fg(n-1-l)$ for $l=0,\dots, n-1$.  Then the decomposition
 \begin{equation}\label{eq:components}
 \fg(n-1-l)=\bigcup_{l(Q)=l} Y_{Q}\cap \fg(n-1-l),
 \end{equation}
is the irreducible component decomposition of the variety $\fg(n-1-l)$, 
where the union is taken over all $K$-orbits $Q$ of length $l$ in $\B$. (cf. Theorem (\ref{eq:theorem})).  
 
 In fact, for $1\leq i\leq j\leq n$ with $j-i=l$, we have
 $$\Ybij\cap\fg(n-1-l)=\Ypij\cap\fg(n-1-l),$$  
 so that  
 \begin{equation}\label{eq:both}
\fg(n-1-l)=\displaystyle\bigcup_{j-i=l} \Ypij\cap\fg(n-1-l).
\end{equation}
  
 \end{thm}
 \begin{proof}
 We first claim that if $l(Q)=l$, then $Y_{Q}\cap\fg(n-1-l)$ is non-empty.
 By Theorem \ref{thm:overlaps}, $Y_{Q}\subset\fg(\geq n-1-l)$. Thus, if $Y_{Q}\cap\fg(n-1-l)$ were empty, then $Y_{Q} \subset \fg(\geq n-l)$.
Hence, by part (1) of Proposition \ref{prop_dimgl}, $\dim(Y_{Q}) \le
n^2-n+l$.   By Equation (\ref{eq:largeYQdim}), $\dim(Y_{Q})=n^2-n+l+1$.
This contradiction verifies the claim.

It follows from Equation (\ref{eq:fgl}) that $\fg(n-1-l)$
 is open in $\fg(\geq n-1-l)$.   Thus, $Y_{Q}\cap\fg(n-1-l)$ is a non-empty Zariski open subset of $Y_{Q}$,
which is irreducible since $Y_{Q}$ is irreducible. 

Now we claim that 
\begin{equation}\label{eq:closeinter}
Y_{Q}\cap\fg(n-1-l)=\overline{Y_{Q}}\cap\fg(n-1-l),
\end{equation}
so that $Y_{Q}\cap\fg(n-1-l)$ is closed 
in $\fg(n-1-l)$.  By Lemma \ref{l:YQclosure},
 $\overline{Y_{Q}}=\bigcup_{Q^{\prime}\subset \overline{Q}}  Y_{Q^{\prime}}$.  
Hence, if (\ref{eq:closeinter}) were not an equality, there
would be $Q^{\prime}$ with $l(Q^{\prime})< l(Q)$ and
$Y_{Q^{\prime}} \cap \fg(n-1-l)$ nonempty.  This contradicts
Theorem \ref{thm:overlaps}, which asserts that $Y_{Q^{\prime}}
\subset \fg(\geq n-l)$, and hence verifies the claim.
It follows that $Y_{Q}\cap\fg(n-1-l)$ is an irreducible, closed subvariety 
of $\fg(n-1-l)$ of dimension $\dim Y_{Q}=\dim\fg(n-1-l)$.
Thus, $Y_{Q}\cap\fg(n-1-l)$ is an irreducible component of $\fg(n-1-l)$.  

Since $l(Q)=l$,  Remark \ref{r:length} implies that
 $Q=Q_{i,j}$ for some $i\leq j$ with $j-i=l$.  Then by Proposition \ref{p:diffYQs} and Equation (\ref{eq:closeinter}),
\begin{equation}\label{eq:YQsequal}
\Ybij\cap\fg(n-1-l)=\Ypij\cap\fg(n-1-l).
\end{equation}
 
Let $Z$ be an irreducible component of $\fg(n-1-l)$.  The proof will
be complete once we show that
 $Z=\Ypij\cap\fg(n-1-l)$ for some $i,j$ with $j-i=l$.
 To do this, 
 consider the nonempty open set
$$U:=\{x\in\fg :\; x_{n-1}\mbox{ is regular semisimple}\}.$$
Let $\widetilde{U}(n-1-l):=\fg(n-1-l)\cap U$.    

Since $\Phi_{n}:\fg\to V^{n-1,n}$ is surjective (by
Remark \ref{rem_kwsurj}), it follows that  $\widetilde{U}(n-1-l)$
is a nonempty Zariski open set of $\fg(n-1-l)$.  By part (2) of Proposition \ref{prop_kwflat}
and Exercise III.9.1 of \cite{Ha}, $\Phi_{n}(U)\subset V^{n-1,n}$ is open. Thus, 
$V^{n-1,n}(n-1-l)\setminus \Phi_{n}(U)$ is a proper, closed subvariety of $V^{n-1,n}(n-1-l)$
and therefore has positive codimension by Lemma \ref{lem:vnirreducible}.  
It follows by part (2) of Proposition \ref{prop_kwflat} and Corollary III.9.6 of \cite{Ha} that 
$\fg(n-1-l)\setminus\widetilde{U}(n-1-l)=\Phi^{-1}_{n}(V^{n-1,n}(n-1-l)\setminus \Phi_{n}(U))$ 
is a proper, closed subvariety of $\fg(n-1-l)$ of positive codimension.  Since $\fg(n-1-l)$ is equidimensional,
it follows that $Z\cap\widetilde{U}(n-1-l)$ is nonempty. 
Thus, it suffices to show that 
\begin{equation}\label{eq:tricky}
\widetilde{U}(n-1-l)\subset\bigcup_{j-i=l} \Ypij\cap\fg(n-1-l).
\end{equation}  

To prove Equation (\ref{eq:tricky}), we consider the following subvariety of $\widetilde{U}(n-1-l)$:
\begin{equation}\label{eq:XI}
\Xi=\{x\in\tilde{U}(n-1-l):\; x_{n-1}=\mbox{diag}[h_{1},\dots, h_{n-1}],\mbox{ and } \sigma(x_{n-1})\cap\sigma(x)=\{h_{1},\dots, h_{n-1-l}\}\} 
\end{equation}
It is easy to check that any element of  $\widetilde{U}(n-1-l)$ is $K$-conjugate
to an element in $\Xi$.  By a linear algebra calculation from Proposition 5.9 of \cite{Col1}, 
elements of $\Xi$ are matrices of the form 

\begin{equation}\label{eq:easymatrix}
\left[\begin{array}{cc}
\begin{array}{cccc}
h_{1}& 0 & \cdots & 0\\
0& h_{2}&\ddots &\vdots\\
\vdots&\;& \ddots & 0\\
0& \cdots &\cdots &h_{n-1}\\
\end{array} & \begin{array}{c}
y_{1}\\
\vdots\\
\vdots\\
y_{n-1}\end{array}\\
\begin{array}{cccc}
z_{1}&\cdots&\cdots&z_{n-1}
\end{array} & w
\end{array}\right], 
\end{equation}
with $h_{i}\neq h_{j}$ for $i\neq j$ and satisfying the equations:

 \begin{equation}\label{eq:nonclosedoverlaps1}
 \begin{split}
 z_{i}y_{i}&=0 \mbox{ for } 1\leq i\leq n-1-l\\
 z_{i}y_{i}&\in \C^{\times} \mbox{ for } n-l\leq i\leq n-1.
 \end{split}
 \end{equation}
 Since the varieties $\Ypij\cap\fg(n-1-l)$ are $K$-stable, it suffices to prove 
\begin{equation}\label{eq:tricky2}
\Xi\subset\bigcup_{j-i=l} \Ypij\cap \fg(n-1-l).
\end{equation}
 
 To prove (\ref{eq:tricky2}), we need to understand the irreducible components of $\Xi$.  For $i=1,\dots, n-1-l$, we define an index $j_{i}$ which takes on two values $j_{i}=U$ ($U$ for upper) or $j_{i}=L$ ($L$ for lower).  Consider the subvariety $\Xi_{j_{1},\dots, j_{n-1-l}}\subset \Xi$ defined by:
 \begin{equation}\label{eq:UorL}
\Xi_{j_{1},\dots, j_{n-1-l}}:=\{ x\in\Xi: z_{i}=0\mbox{ if } j_{i}=U, y_{i}=0 \mbox{ if } j_{i}=L\}.
\end{equation}
    Then
      \begin{equation}\label{eq:upsunion}
  \Xi= \bigcup_{j_{i}=U,\, L} \Xi_{j_{1},\dots, j_{n-1-l}}.
  \end{equation} 
  is the irreducible component decomposition of $\Xi$.  
  
  We now consider the irreducible variety $\Xi_{j_{1},\dots, j_{n-1-l}}$.  Suppose that for the subsequence $1\leq i_{1}<\dots<i_{k-1}\leq n-1-l$ we have $j_{i_{1}}=j_{i_{2}}=\dots=j_{i_{k-1}}=U$ and that for the complementary subsequence $i_{k}<\dots< i_{n-1-l}$ we have $j_{i_{k}}=j_{i_{k+1}}=\dots=j_{i_{n-1-l}}=L$.  Then a simple computation with flags shows that elements of the variety $\Xi_{j_{1},\dots, j_{n-1-l}}$ stabilize the $n-l$-step partial flag in $\C^{n}$
   \begin{equation}\label{eq:partial}
  e_{i_{1}}\subset e_{i_{2}}\subset \dots\subset e_{i_{k-1}}\subset\underbrace{ e_{n-l},\dots, e_{n-1}, e_{n}}_{k}\subset e_{i_{k}}\subset e_{i_{k+1}}\subset \dots\subset e_{i_{n-1-l}}.
     \end{equation}
     (If $l=0$ the partial flag in (\ref{eq:partial}) is a full flag with $e_{n}$ in the $k$-th position.)
     It is easy to see that there is an element of $K$ that maps the partial flag in Equation (\ref{eq:partial}) 
     to the partial flag $\mathcal{P}_{k,k+l}$ in Equation (\ref{eq:firstpartial}):
\begin{equation}\label{eq:secondpartial}
 \mathcal{P}_{k,k+l}=(e_{1}\subset e_{2}\subset\dots\subset e_{k-1}\subset\underbrace{e_{k},\dots, e_{k+l-1}, e_{n}}_{k}\subset e_{k+l}\subset\dots\subset e_{n-1}).
 \end{equation}
 (If $l=0$ the partial flag $\mathcal{P}_{k,k+l}$ is the full flag $\mathcal{F}_{k,k}$ (see Equation (\ref{eq:flagi})).)
  Thus, $\Xi_{j_{1},\dots, j_{n-1-l}}\subset Y_{{\fp_{k,k+l}}}\cap \fg(n-1-l)$.  Equation (\ref{eq:upsunion}) then implies that $\Xi\subset \bigcup_{j-i=l} \Ypij\cap\fg(n-1-l)$.  
 \end{proof}
 
Using Theorem \ref{thm:king}, we can obtain the irreducible component decomposition of the variety $\fg(\geq n-1-l)$ for
any $l=0,\dots, n-1$.  
\begin{cor}\label{c:geq}
The irreducible component decomposition of the variety $\fg(\geq n-1-l)$ is
\begin{equation}\label{eq:geq}
\fg(\geq n-1-l)=\bigcup_{j-i=l} \Ypij=\bigcup_{l(Q)=l} \overline{Y_{Q}}.
\end{equation}
\end{cor}
\begin{proof}
Taking Zariski closures in Equation (\ref{eq:both}), we obtain
\begin{equation}\label{eq:closure}
\overline{\fg(n-1-l)}=\bigcup_{j-i=l} \overline{\Ypij\cap\fg(n-1-l)}
\end{equation}
is the irreducible component decomposition of the variety $\overline{\fg(n-1-l)}$.  
By Proposition \ref{prop_dimgl}, $\overline{\fg(n-1-l)}=\fg(\geq n-1-l)$, 
and  by Theorem \ref{thm:overlaps}, $\Ypij\subset \fg(\geq n-1-l)$. 
Hence $\Ypij\cap \fg(n-1-l)$ is Zariski open in the irreducible
variety $\Ypij$, and is nonempty by Theorem \ref{thm:king}. Therefore $\overline{\Ypij\cap\fg(n-1-l)}=\Ypij$.  Equation (\ref{eq:geq})
now follows from Equation (\ref{eq:closure}) and Proposition \ref{p:diffYQs}.  
\end{proof}

Theorem \ref{thm:king} says something of particular interest to linear algebraists 
in the case where $l=0$.  It states that the variety $\fg(n-1)$ consisting of elements $x\in\fg$ where
the number of coincidences in the spectrum between $x_{n-1}$ and $x$ is maximal can be described in terms of 
closed $K$-orbits on $\mathcal{B}$, which are the $K$-orbits $Q$ with
$l(Q)=0$. 
 It thus connects the most degenerate case of spectral coincidences to the 
simplest $K$-orbits on $\mathcal{B}$.  More precisely, we have:

\begin{cor}\label{cor:degen}
The irreducible component decomposition of the variety $\fg(n-1)$ is 
$$
\fg(n-1)=\bigcup_{l(Q)=0} Y_{Q}. 
$$
\end{cor}

Using Corollary \ref{cor:degen} and Theorem \ref{thm:completeinter}, we obtain a precise description of the irreducible components 
of the variety $SN_{n}$ introduced in Equation (\ref{eq:strongnil}). 
\begin{prop} 
\label{prop:irredsnn}
Let $\fb_{i,i}$ be the Borel subalgebra of $\fg$ which stabilizes the flag $\mathcal{F}_{i,i}$ 
in Equation (\ref{eq:flagi}) and let $\fn_{i,i}=[\fb_{i,i},\fb_{i,i}]$.  
The irreducible component decomposition of $SN_{n}$ is given by:
\begin{equation}\label{eq:cptsSNn}
SN_{n}=\bigcup_{i=1}^{n} \Ad(K)\fn_{i,i},
\end{equation}
where $\Ad(K)\fn_{i,i}\subset\fg$ denotes the $K$-saturation of $\fn_{i,i}$ in $\fg$.

\end{prop}
\begin{proof}
We first show that $\Ad(K) \fn_{i,i}$ is an irreducible component 
of $SN_{n}$ for $i=1,\dots, n$.  A simple computation using the flag $\mathcal{F}_{i,i}$ in 
Equation (\ref{eq:flagi}) shows that $\fn_{i,i}\subset SN_{n}$.  Since $SN_{n}$ is $K$-stable, it 
follows that $\Ad(K) \fn_{i,i}\subset SN_{n}$.  

Recall the Grothendieck resolution $\widetilde{\fg}=\{(x,\fb): x\in\fb\}\subset\fg\times\B$ and 
the morphisms $\pi:\widetilde{\fg}\to \B$, $\pi(x,\fb)=\fb$ and $\mu:\widetilde{\fg}\to \fg$, $\mu(x,\fb)=x$. 
  Let $Q_{i,i}=K\cdot\fb_{i,i}\subset\B$ be the $K$-orbit 
through $\fb_{i,i}$.  Corollary 3.1.33 of \cite{CG} gives a $G$-equivariant isomorphism 
$\widetilde{\fg}\cong G\times_{B_{i,i}}\fb_{i,i}$.  Under this isomorphism 
$\pi^{-1}(Q_{i,i})$ is identified 
with the closed subvariety $K\times_{K\cap B_{i,i}}\fb_{i,i}\subset G\times_{B_{i,i}}\fb_{i,i}$.  
The closed subvariety $K\times_{K\cap B_{i,i}}\fn_{i,i}\subset K\times_{K\cap B_{i,i}}\fb_{i,i}$ maps surjectively under $\mu$
to $\Ad(K) \fn_{i,i}$.  Since $\mu$ is proper, $\Ad(K) \fn_{i,i}$ is closed and irreducible.  
We also note that the restriction of $\mu$ to $K\times_{K\cap B_{i,i}}\fn_{i,i}$ generically has finite fibers 
(Proposition 3.2.14 of \cite{CG}).  Thus, the same reasoning that we used in Equation (\ref{eq:largeYQdim})  shows that 
\begin{equation}\label{eq:nulldim}
\dim \Ad(K) \fn_{i,i}=\dim(K\times_{K\cap B_{i,i}}\fn_{i,i})=\dim(Y_{Q_{i,i}})-\mbox{rk}(\fg)=d_{n},
\end{equation}
where $\mbox{rk}(\fg)$ denotes the rank of $\fg$.
Thus, by Theorem \ref{thm:completeinter}, $\Ad(K) \fn_{i,i}$ is an irreducible component of 
$SN_{n}$.  

We now show that every irreducible component of $SN_{n}$ is of the form $\Ad(K) \fn_{i,i}$ for some $i=1,\dots, n$. 
It follows from definitions
that $SN_{n}\subset \fg(n-1)\cap \mathcal{N}$, where $\mathcal{N}\subset\fg$ is the nilpotent cone in $\fg$.  
Thus, if $\X$ is an irreducible component of $SN_{n}$, 
then $\X\subset \Ad(K) \fn_{i,i}$ by Corollary \ref{cor:degen}.  But then 
$\X=\Ad(K) \fn_{i,i}$ by Equation (\ref{eq:nulldim}) and Theorem \ref{thm:completeinter}.

\end{proof}

We say that an element $x\in \fg$ is {\it n-strongly regular} if
the set $$dJZ_{n}(x):=\{ df_{i,j}(x) : i=n-1, n; j=1, \dots, i \}$$ is linearly independent
in the cotangent space $T_{x}^{*}(\fg)$ of $\fg$ at $x$.  We view $\fg_{n-1}$ as the top lefthand corner of $\fg$.  It follows from a well-known result of Kostant (Theorem 9 of \cite{Ko}) that $x_{i}\in\fg_{i}$ is regular if and only if the set $\{ df_{i,j}(x) :  j=1, \dots, i \}$ is linearly independent.   
If $x_{i}\in\fg_{i}$ is regular, and we identify $T_{x}^{*}(\fg)=\fg^{*}$ with $\fg$ using the trace form $<<x,y>>=tr(xy)$, then 
$$
\mbox{ span }\{ df_{i,j}(x):j=1, \dots, i\}=\fz_{\fg_{i}}(x_{i}),
$$
where $\fz_{\fg_{i}}(x_{i})$ denotes the centralizer of $x_{i}$ in $\fg_{i}$.  Thus, it follows that
 $x\in\fg$ is $n$-strongly regular if and only if $x$ satisfies the following two conditions:
\begin{equation}\label{eq:nsr}
\begin{split}
&(1)\; x\in \fg \mbox { and }  x_{n-1}\in\fg_{n-1}\mbox{ are regular; and }\\
&(2)\; \fz_{\fg_{n-1}}(x_{n-1})\cap\fz_{\fg}(x)=0. \end{split}
\end{equation}
\begin{rem}\label{r:notradical}
We claim that the ideal $I_{n}$ is radical if and only if $n\le 2$.
The assertion is clear for $n=1$, and we assume $n\ge 2$ in the sequel.
Indeed, by Theorem 18.15(a) of \cite{Eis}, the ideal $I_{n}$ is radical
if and only if the set $dJZ_{n}$ is linearly independent on a dense open
set of  each irreducible
component of $SN_{n}=V(I_{n})$.  It follows that $I_{n}$ is radical if and
only if each irreducible component of $SN_{n}$ contains n-strongly
regular elements.  Let $\fn_{+}=[\fb_{+}, \fb_{+}]$ and $\fn_{-}=[\fb_{-},\fb_{-}]$
be the strictly upper and lower triangular matrices, respectively.  
By Proposition \ref{prop:irredsnn} above, $SN_{n}$ has exactly $n$
irreducible components.  It follows from the discussion after Equation (\ref{eq:flagi}) that two 
of them are
$K\cdot \fn_{+}$ and $K\cdot \fn_{-}$.  
By Proposition 3.10 of \cite{CEKorbs}, the only
irreducible components of $SN_{n}$ which contain n-strongly regular
elements are $K\cdot \fn_{+}$ and $K\cdot \fn_{-}$.  The claim now
follows.  See Remark 1.1 of \cite{Ov} for a related observation, which
follows also from the analysis proving our claim.
\end{rem}


 
\bibliographystyle{amsalpha.bst}

\bibliography{bibliography-1}

\def\cprime{$'$} \def\cprime{$'$} \def\cprime{$'$} \def\cprime{$'$}
  \def\cprime{$'$} \def\cprime{$'$} \def\cprime{$'$} \def\cprime{$'$}
  \def\cprime{$'$} \def\cprime{$'$} \def\cprime{$'$} \def\cprime{$'$}
  \def\cprime{$'$} \def\cprime{$'$}
\providecommand{\bysame}{\leavevmode\hbox to3em{\hrulefill}\thinspace}
\begin{thebibliography}{Hum75}

\bibitem[BH00]{BH}
Michel Brion and Aloysius~G. Helminck, \emph{On orbit closures of symmetric
  subgroups in flag varieties}, Canad. J. Math. \textbf{52} (2000), no.~2,
  265--292.

\bibitem[Bri87]{Brionclass}
M.~Brion, \emph{Classification des espaces homog\`enes sph\'eriques},
  Compositio Math. \textbf{63} (1987), no.~2, 189--208.

\bibitem[CE]{CEexp}
Mark Colarusso and Sam Evens, \emph{The {G}elfand-{Z}eitlin integrable system
  and {K}-orbits on the flag variety}, to appear in: ``Symmetry: Representation
  Theory and its Applications: In Honor of Nolan R. Wallach," Progr. Math.
  Birkauser, Boston.

\bibitem[CE12]{CEKorbs}
Mark Colarusso and Sam Evens, \emph{K-orbits on the flag variety and strongly
  regular nilpotent matrices}, Selecta Math. (N.S.) \textbf{18} (2012), no.~1,
  159--177.

\bibitem[CG97]{CG}
Neil Chriss and Victor Ginzburg, \emph{Representation theory and complex
  geometry}, Birkh\"auser Boston Inc., Boston, MA, 1997.

\bibitem[Col11]{Col1}
Mark Colarusso, \emph{The orbit structure of the {G}elfand-{Z}eitlin group on
  {$n\times n$} matrices}, Pacific J. Math. \textbf{250} (2011), no.~1,
  109--138.

\bibitem[Eis95]{Eis}
David Eisenbud, \emph{Commutative algebra}, Graduate Texts in Mathematics, vol.
  150, Springer-Verlag, New York, 1995, With a view toward algebraic geometry.

\bibitem[FO05]{FuOvsfa}
Vyacheslav Futorny and Serge Ovsienko, \emph{Kostant's theorem for special
  filtered algebras}, Bull. London Math. Soc. \textbf{37} (2005), no.~2,
  187--199.

\bibitem[Gro03]{SGA1}
Alexander Grothendieck, \emph{Rev\^{e}tements \'{E}tales et groupe fondamental
  ({SGA} 1)}, Documents Math\'ematiques (Paris) [Mathematical Documents
  (Paris)], 3, Soci\'et\'e Math\'ematique de France, Paris, 2003, S{\'e}minaire
  de G{\'e}om{\'e}trie Alg{\'e}brique du Bois Marie, 1960-1961, Augment{\'e} de
  deux expos{\'e}s de Mich{\`e}le Raynaud. [With two expos{\'e}s by Mich{\`e}le
  Raynaud].

\bibitem[Har77]{Ha}
Robin Hartshorne, \emph{Algebraic geometry}, Springer-Verlag, New York, 1977,
  Graduate Texts in Mathematics, No. 52.

\bibitem[Has04]{Hashi}
Takashi Hashimoto, \emph{{$B_{n-1}$}-orbits on the flag variety {${\rm
  GL}_n/B_n$}}, Geom. Dedicata \textbf{105} (2004), 13--27.

\bibitem[Hum75]{Hum}
James~E. Humphreys, \emph{Linear algebraic groups}, Springer-Verlag, New York,
  1975, Graduate Texts in Mathematics, No. 21.

\bibitem[Kos63]{Ko}
Bertram Kostant, \emph{Lie group representations on polynomial rings}, Amer. J.
  Math. \textbf{85} (1963), 327--404.

\bibitem[KW06]{KW1}
Bertram Kostant and Nolan Wallach, \emph{Gelfand-{Z}eitlin theory from the
  perspective of classical mechanics. {I}}, Studies in {L}ie theory, Progr.
  Math., vol. 243, Birkh\"auser Boston, Boston, MA, 2006, pp.~319--364.

\bibitem[Mat79]{Mat79}
Toshihiko Matsuki, \emph{The orbits of affine symmetric spaces under the action
  of minimal parabolic subgroups}, J. Math. Soc. Japan \textbf{31} (1979),
  no.~2, 331--357.

\bibitem[Mat86]{Mat}
Hideyuki Matsumura, \emph{Commutative ring theory}, Cambridge Studies in
  Advanced Mathematics, vol.~8, Cambridge University Press, Cambridge, 1986,
  Translated from the Japanese by M. Reid.

\bibitem[M{\=O}90]{MO}
Toshihiko Matsuki and Toshio {\=O}shima, \emph{Embeddings of discrete series
  into principal series}, The orbit method in representation theory
  ({C}openhagen, 1988), Progr. Math., vol.~82, Birkh\"auser Boston, Boston, MA,
  1990, pp.~147--175.

\bibitem[Ovs03]{Ov}
Serge Ovsienko, \emph{Strongly nilpotent matrices and {G}elfand-{Z}etlin
  modules}, Linear Algebra Appl. \textbf{365} (2003), 349--367, Special issue
  on linear algebra methods in representation theory.

\bibitem[RS90]{RS}
R.~W. Richardson and T.~A. Springer, \emph{The {B}ruhat order on symmetric
  varieties}, Geom. Dedicata \textbf{35} (1990), no.~1-3, 389--436.

\bibitem[Ser00]{Serrelocalg}
Jean-Pierre Serre, \emph{Local algebra}, Springer Monographs in Mathematics,
  Springer-Verlag, Berlin, 2000, Translated from the French by CheeWhye Chin
  and revised by the author.

\bibitem[VK78]{KimelVin}
{\'E}.~A. Vinberg and B.~N. Kimel{\cprime}fel{\cprime}d, \emph{Homogeneous
  domains on flag manifolds and spherical subsets of semisimple {L}ie groups},
  Funktsional. Anal. i Prilozhen. \textbf{12} (1978), no.~3, 12--19, 96.

\bibitem[Yam97]{Yam}
Atsuko Yamamoto, \emph{Orbits in the flag variety and images of the moment map
  for classical groups. {I}}, Represent. Theory \textbf{1} (1997), 329--404
  (electronic).

\end{thebibliography}


\end{document}